\documentclass[11pt,reqno,tbtags]{amsart}
\usepackage{amssymb,epsf,epsfig}
\numberwithin{equation}{section}

\makeatletter
\newtheorem*{rep@theorem}{\rep@title}
\newcommand{\newreptheorem}[2]{%
\newenvironment{rep#1}[1]{%
 \def\rep@title{#2 \ref{##1}}%
 \begin{rep@theorem}}%
 {\end{rep@theorem}}}
\makeatother

\newtheorem{theorem}{Theorem}[section]
\newreptheorem{theorem}{Theorem}

\theoremstyle{definition}

\theoremstyle{remark}

\newcounter{thmenumerate}

\newcounter{xenumerate}

\newcommand\E{\operatorname{\mathbb E{}}}

\renewcommand\Pr{\operatorname{\mathbb P{}}}

\newcommand\eps{\varepsilon}

\renewcommand\phi{\varphi}

\newcommand\cB{\mathcal B}

\newcommand\cF{\mathcal F}
\newcommand\cG{\mathcal G}

\newcommand\N{{\mathbb N}}

\newcommand\R{{\mathbb R}}
\newcommand\M{{\mathbb M}}

\begin{document}
\title{The Mathematics of Causal Sets}

\date{\today} 

\author{Graham Brightwell}
\address{Department of Mathematics, London School of Economics,
Houghton Street, London WC2A 2AE, United Kingdom}
\email{g.r.brightwell@lse.ac.uk}
\urladdr{http://www.maths.lse.ac.uk/Personal/graham/}

\author{Malwina Luczak} \thanks{The research of Malwina Luczak is supported by an EPSRC Leadership Fellowship, grant reference EP/J004022/2.}
\address{School of Mathematics, Queen Mary, University of London}
\email{m.luczak@qmul.ac.uk}
\urladdr{http://www.maths.qmul.ac.uk/$\sim$luczak/}


\maketitle

\begin{abstract}
The causal set approach to quantum gravity is based on the hypothesis that the underlying structure of spacetime is that of a random partial order.
We survey some of the interesting mathematics that has arisen in connection with the causal set hypothesis, and describe how the mathematical
theory can be translated to the application area.  We highlight a number of open problems of interest to those working in causal set theory.
\end{abstract}

\section {The causal set hypothesis}

This article constitutes a survey of the mathematical results related to the causal set approach to quantum gravity.  Simply put, the
causal set hypothesis is that spacetime is fundamentally a discrete structure, consisting of a set of points equipped with a partial order.  On a large
scale, such a discrete structure may approximate a Lorentzian spacetime manifold, with the partial order given by ``is in the past light cone of'',
as in the theory of general relativity.
On a small scale, the behaviour will be completely different, with an inherent discreteness and ``randomness'', as in the theory of quantum mechanics.
Causal set theory is one of a number of proposed alternatives to string theory; research took off in the late 1980s following work of
Bombelli, Lee, Meyer and Sorkin~\cite{BLMS}, but the concept occurred to others independently, notably in the preprint of
Myrheim~\cite{Myrheim} from 1978.  Although causal set theory is now well-established as a field, it does not
enjoy the same status or support as string theory, and progress has accordingly been less fast.

Underlying the theory is the notion that the ``correct'' model of spacetime is a model of a random partially ordered set.  For each of the most common
models of random partial orders in the mathematical literature -- see Brightwell~\cite{Bsurvey} for a survey -- some variant or extension has its
attractions as a candidate for a model of spacetime.  Our aim in this article is to review the mathematics, explain its relevance to physics, and discuss
some potentially important open problems.

There are two main strands to the physics.  One is to argue that we know what the discrete set should look like: it should be a set of points embedded
``faithfully'' (i.e., with about the same density in almost all regions of equal volume) in a 4-dimensional Lorentzian manifold, which on a moderate
scale is essentially flat.  This can be achieved, at least approximately, by taking the points of our discrete set to be a Poisson process in
4-dimensional Minkowski space $\M^4$, with the order induced on this set from the underlying manifold.  We then ask what we can do with such a discrete
partially ordered set: can we in any sense recover the structure of the manifold it came from?

The other strand aims to address the question of ``why'' the spacetime resembles $\M^4$ on large scales.  The ultimate prize would be a simple, natural
random process, making no use of any underlying manifold, whose samples do resemble discrete sets that are ``faithfully'' embedded in $\M^4$.  There is
currently no candidate for such a process, but much interesting mathematics has resulted from various approaches to the issue.

Given that spacetime is to be discrete, one might still ask why it makes sense to use a random structure?  The most obvious alternative is to model space as a lattice-like structure, with layers representing time.  However, a model of this form will exhibit ``preferred directions'' along the lattice that we do not
observe in practice.  In particular, observations suggest that the spacetime structure of the universe is invariant under Lorentz transformations, and so
our discrete model ought also to enjoy this property, or a very close approximation of it.  More generally, there is no explicit construction known of a
discrete subset of $\M^4$ that has approximately the ``right'' number of points in each interval.  As we shall discuss later, taking a Poisson process in $\M^4$
does not completely solve this problem, but at least it provides the basis for a theory that does not discriminate between intervals of equal volume.
There is also the intuition coming from quantum mechanics that aspects of reality are best modelled using some ``randomness'', albeit not necessarily
one involving a classical probability measure.

There are a number of surveys of the theory of causal sets from a physics perspective~\cite{Sorkin-survey, Henson-survey, Surya-survey, Dowker-survey},
and the reader is referred to those for more background and motivation.  The present article is not intended as a full survey
of the causal set programme; two particular omissions are (i)~the prediction from the theory of a small but non-zero value of the cosmological constant,
subsequently supported by observations, and (ii)~recent progress on the important project of developing notions of fields and particles on top of the
underlying causal set.

A partially ordered set, or {\em poset}, is a pair $(X,<)$ where $<$ is a transitive, irreflexive relation on the {\em ground-set} $X$.
A causal set is intended to model a universe with a big bang, meaning that the past light cone of every point has finite measure.  The equivalent statement
in the discrete setting is that every element $x$ has only finitely many elements below it.  It is customary in the physics literature to take this as the
definition of a {\em causal set} $(X,<)$: a poset with countable ground-set $X$ such that, for every $x \in X$, the set $\{ y \in P : y \le x\}$
is finite.  For many purposes, one might instead work with a {\em locally finite} poset, where the {\em intervals} $[x,y] = \{ z : x \le z \le y \}$
are finite.  We introduce other definitions as we proceed.

Whether or not causal set theory is ultimately a success, even the possibility of an application area enriches the study of random partial orders, and provides
an interesting alternative viewpoint regarding what questions are of natural interest.

\section{Poisson processes in Lorentzian manifolds} \label{sec:two}

In this section, we examine models arising by taking a continuum structure equipped with a partial order and a compatible measure, generating a discrete set
of points by taking a Poisson process in the continuum, and then looking purely at the discrete set with its induced partial order.

The most basic continuum model of spacetime is 4-dimensional Minkowski space $\M^4$.  Minkowski space $\M^d$ is a Lorentzian manifold, with $d-1$ dimensions
of space and one of time: one can represent $\M^d$ as $\R^{d-1} \times \R$, with the product Lebesgue measure, and an order relation given by
$({\bf x},s) \le ({\bf y},t)$ if $|{\bf x} - {\bf y}| \le t-s$.  Here we have expressed our coordinates in units where the speed of light is~1, so to
say that one point is below another means that light may travel from the lower point of spacetime to the higher.  The set of points above
${\bf a} = ({\bf x},s)$ is the {\em future light cone} of ${\bf a}$, and the set of points below ${\bf a}$ is its {\em past light cone}.
A pair of points is said to be {\em timelike} if they are comparable in the partial order:
in this case the distance between $({\bf x},s)$ and $({\bf y},t)$ is given by $\sqrt{(t-s)^2 - |{\bf x} - {\bf y}|^2}$ and is called the {\em proper time}
between the two points.  Note that this distance might
be zero: this corresponds to one point lying on the boundary of the future light cone of the other.  If a pair of points is incomparable, they are said to be
{\em spacelike}, and the distance between them is defined as $\sqrt{|{\bf x} - {\bf y}|^2 - (t-s)^2}$.
The measure, the partial order, and the distance function on this space are then invariant under the group of
{\em Lorentz transformations}.  These include isometries of the space $\R^{d-1}$, with the time axis left fixed, and also the {\em Lorentz boosts}:
with respect to any orthogonal coordinate system for the space $\R^{d-1}$, a boost in the direction of the first coordinate is given by
mapping $((x_1, \dots, x_{d-1}),t)$ to $((\gamma(x_1 -\beta t),x_2, \dots, x_{d-1}),\gamma(t - \beta x_1))$, where $\beta \in (-1,1)$ is the
{\em velocity coefficient} (measuring the magnitude of the boost) and $\gamma = (1-\beta^2)^{-1/2}$ is the {\em Lorentz factor}.  A boost represents
transforming the time axis of one observer to that of another travelling at speed $\beta$ in the direction of the first coordinate with respect to the
first observer.

Although it is natural to consider 4-dimensional models of spacetime, it is also reasonable to think about simpler ``toy models'' where, for instance,
$\M^2$ replaces $\M^4$.  It is also desirable to be able to replace $\M^4$, representing a flat spacetime, with a more
general Lorentzian manifold with curvature.

There is an alternative visual interpretation of a finite partial order $P$ embedded in $\M^d$.
Imagine a ``horizontal'' plane $H$ to the past of all points of $P$ and, for each point $x$ of $P$, consider the intersection $S_x$ of the past
light cone of $x$ with $H$.  This intersection will be a $(d-1)$-sphere, and we have that $S_x \subseteq S_y$ if and only if $x$ is below $y$ in the order
on $\M^d$.  Hence the collection of spheres $S_x$, ordered by containment, is a representation of the partial order $P$.  This construction can be
reversed, so that a finite partial order can be embedded in $\M^d$ if and only if it has a representation as a collection of $(d-1)$-spheres ordered
by containment.  A partial order with such a representation is called a {\em circle order} if $d=3$, and a {\em $(d-1)$-sphere order} for higher~$d$.
See, for instance, Meyer~\cite{Meyer}, or Brightwell and Winkler~\cite{BW-sphere}.

For two points ${\bf a}, {\bf b}$ in Minkowski space with ${\bf a} \le {\bf b}$, the interval $[{\bf a},{\bf b}]$ is given by
$\{ {\bf z} : {\bf a} \le {\bf z} \le {\bf b} \}$.  For two intervals of the same positive measure, there is a Lorentz transformation taking one to the other.
Thus any two intervals of positive measure are ``equivalent'' up to a scale factor; for most purposes, we can restrict our attention to one such
interval.

Any (discrete) subset of $\M^d$ inherits the order structure from the Minkowski space, and we can then ask questions about the partial orders arising.
There are various results concerning how much of the structure of a Lorentzian manifold can be recovered from just the causal (partial order) structure,
and the basic conclusion is that the topology and geometry can be recovered completely, up to a conformal factor specifying the volume measure.
The most general results in this line are those of Hawking, King and McCarthy~\cite{Hawking} and Malament~\cite{Malament}.

We now take a Poisson process $P$ with intensity $n$ in $\M^d$.  The partial order on $\M^d$ then induces a partial order on the discrete set of points, and we
can look at this purely as a discrete structure, ``rubbing out'' the manifold used to construct it.  The basic principle is that the order relation from
the manifold survives in the discrete approximation, and the volume measure can be approximated by counting the points in any given region.  Hence one
might hope that the discrete structure contains enough of the structure of the original manifold to form a sufficient replacement for it.

From a mathematical perspective, one may fix some interval $[{\bf a}, {\bf b}]$ in the manifold, let the intensity $n$ of the Poisson process tend to
infinity, and look at the limiting distribution of suitably scaled parameters of the partial order $P[{\bf a},{\bf b}]$ restricted to the interval.
The physical perspective does not allow us the luxury of taking the limit: physical considerations lead to a proposition for
a natural value of the intensity, where the expected number of points in a ``standard'' volume of spacetime, with spatial volume on the order of a cubic
metre and time duration on the order of the time for light to travel one metre, is about $10^{136}$: see~\cite{Sorkin-survey} for a detailed
argument.  Indeed, the point of the theory is that small intervals of spacetime should be unlike the continuum, so that we can
see ``quantum behaviour'' and ``relativistic behaviour'' in the same model, but at different scales.

A fundamental question is whether, if we go to a large enough scale, parameters of the partial order are good approximations to parameters of the underlying
manifold.  Versions of this question were framed in the seminal paper of Bombelli, Lee, Meyer and Sorkin~\cite{BLMS}, as well as in the earlier preprint
of Myrheim~\cite{Myrheim}.

Consider two points ${\bf a}$ and ${\bf b}$ in the Poisson process on the underlying manifold $M$.  They will be at some distance in $M$, which can
be either timelike, if ${\bf a} \le {\bf b}$, or spacelike, if the two points are incomparable in the partial order on $M$.  A natural candidate for an
approximation of the (scaled) distance between two timelike points is the {\em height} of the partial order $P[{\bf a},{\bf b}]$, i.e., the length~$h$
of the longest chain ${\bf a} = {\bf x}_0 < {\bf x}_1 < \cdots < {\bf x}_h = {\bf b}$.  The idea that the height could serve as a discrete version of timelike
distance was first proposed by Myrheim~\cite{Myrheim}.
More generally, one can ask about the distribution of the height of the random partial order $P[{\bf a}, {\bf b}]$.

There is a closely related model of random partial orders that has received much more attention in the mathematical literature.  Instead of Minkowski space,
consider $\R^d$, with Lebesgue measure and co-ordinate order, so $(a_1, a_2, \dots, a_d) \le (b_1, b_2, \dots, b_d)$ if and only
if $a_i \le b_i$ for each~$i$.  Again we consider a Poisson process of intensity $n$ in $\R^d$, and take the induced partial order on the discrete set of
points of the Poisson process.  For many purposes, it is convenient to work with the finite unit cube $[0,1]^d$, i.e., the interval 
$[(0,\dots,0),(1,\dots,1)]$ instead.

For $d=2$, the map $f: \M^2 \to\R^2$ given by $f(x,t) = \frac{1}{\sqrt 2} (t+x,t-x)$ is an isomorphism, but, for higher~$d$, $\M^d$ and $\R^d$ are 
different spaces, and the random partial orders have different distributions.
In $\R^d$ with the coordinate order, the ``future light cone'' of $x$, the set of points above $x$ in $\R^d$, has cross-sections that are simplices.
For instance, in $\R^3$, the model has what have been dubbed ``triangular light cones'', with three distinguished directions.

There is an alternative interpretation of the above model.  First, there is little essential difference between taking a Poisson process of intensity~$n$ in
the unit cube, and taking exactly $n$ (labelled) points uniformly and independently in the cube.  The latter amounts to taking $n$ independent $d$-tuples,
with each of the $nd$ numbers uniform in $[0,1]$.  All we are interested in is the partial order on the points, which is the intersection of the
$d$ linear orders given by the coordinates, and these coordinate orders are independent, each being uniform over all $n!$ total orders.  Therefore an
equivalent description is given by taking $d$ linear orders uniformly and independently on an $n$-element ground-set, and intersecting the orders.
The {\em dimension} of a poset $P$ is the minimum number of linear orders on its ground-set whose intersection is $P$, or equivalently the minimum~$d$
such that $P$ can be embedded in $\R^d$.  A random partial order $P_d(n)$ obtained by intersecting $d$ iid uniform linear orders on an $n$-element
ground-set is a {\em random $d$-dimensional partial order}.  This model was introduced by Winkler~\cite{Winkler} in 1985; as we have just seen, essentially
the same model is obtained by taking a Poisson process of intensity~$n$ in $[0,1]^d$.

There is a long history to the study of the height $H_2(n)$ of a random 2-dimensional partial order $P_2(n)$, partly because this has a further guise as the
length of a longest increasing subsequence in a uniformly random permutation of $n$ points.  A number of authors contributed to showing that the height
is approximately $2 \sqrt n$, and then Baik, Deift and Johansson~\cite{BDJ} gave the following impressively precise result.

\begin{theorem} \label{thm:bdj}
$$
\Pr \Big( \frac{H_2(n) - 2n^{1/2}}{n^{1/6}} \le x \Big) \to F(x), \quad \mbox{ as } n\to \infty,
$$
where $F(x)$ is a specified function.  In particular, the mean and standard deviation of $H_2(n)$ satisfy
$$
\E H_2(n) = 2n^{1/2} - \alpha n^{1/6} + o(n^{1/6}); \quad \sigma (H_2(n)) = \beta n^{1/6} + o(n^{1/6}),
$$
where $\alpha \simeq 1.711$, $\beta \simeq 0.902$.
\end{theorem}

See Aldous and Diaconis~\cite{AD} for more on the history and context of this result.

The proper time between two timelike points ${\bf a}$ and ${\bf b}$ in $\M^2$ is, up to a constant factor, $({\rm vol} [{\bf a},{\bf b}])^{1/2}$, and
so the result above can be interpreted as saying that, for a Poisson process of unit intensity in $\M^2$, the height of the partial order
in $[{\bf a}, {\bf b}]$ gives, up to a constant factor, a very good approximation to the proper time between the two points.

The degree of precision in Theorem~\ref{thm:bdj} remains far out of reach for the height of a random partial order in more than~2 dimensions, either in
the cube or in Minkowski space.  Bollob\'as and Winkler~\cite{BolW} showed that the height $H_d(n)$ of $P_d(n)$ behaves as $n^{1/d}$: they showed that
there is a sequence of constants $(c_d)$, with $c_d < e$ for all $d$ and $c_d \to e$ as $d \to \infty$, such that, for each fixed $d$,
$H_d(n)n^{-1/d} \to c_d$ in probability.  Even the values of $c_d$ for $d \ge 3$ remain unknown, let alone any finer estimates.

This result was generalised by Bollob\'as and Brightwell~\cite{BB-scale} to a setting where, as in both $[0,1]^d$ and $\M^d$,
all intervals ``look the same''.  A {\em partially ordered measure space}
is a measure space $(X, \cF, \mu)$ equipped with a partial order $<$ so that all intervals $[x,y]$ are measurable.  Two such spaces 
$(X,\cF,\mu,<)$ and $(X',\cF',\mu',<')$ 
are said to be {\em order-isomorphic} if there is a bijection $\lambda: X \to X'$ respecting the partial order and the $\sigma$-field;
they are said to be {\em scale-isomorphic}
if additionally there is a constant $\alpha$ such that $\mu'(\lambda(A)) = \alpha \mu(A)$ for all $A \in \cF$.  A partially ordered measure space is
{\em homogeneous} if any two intervals of positive measure are scale-isomorphic.  A {\em box space} is a homogeneous partially ordered measure space
that is itself an interval of measure~1.

Both $[0,1]^d$ and an interval in $\M^d$ of measure~1 are box spaces, and Bollob\'as and Brightwell~\cite{BB-scale} give a few other examples as well.
They proved that every box space has a ``dimension'' $d \ge 1$ (which may be infinite, or non-integer) so that, for any interval $[{\bf a},{\bf b}]$ of
positive measure, the maximum, over ${\bf c} \in [{\bf a},{\bf b}]$,
of the minimum of $\mu[{\bf a},{\bf c}]$ and $\mu[{\bf c},{\bf b}]$ is $2^{-d}\mu[{\bf a},{\bf b}]$.

For a Poisson process of density $n$ in a box space $X$, let $H_X(n)$ denote the height of the induced partial order.
The following result of Bollob\'as and Brightwell~\cite{BB-scale} states that the height is asymptotic to $m_X n^{1/d}$ for some constant $m_X$, depending
on the space.

\begin{theorem} \label{thm:BB-scale}
Let $X$ be a box space of dimension $d$.  Then there is a constant $m_X$ such that
$n^{-1/d} H_X(n) \to m_X$ in probability.
\end{theorem}

The only case of interest to us where $m_X$ is known is $X = [0,1]^2$, where $m_X = 2$.
For $X$ an interval in $\M^d$, $d \ge 3$, the bounds given in~\cite{BB-scale} are
$$
\frac{2^{1-1/d}}{\Gamma(1+1/d)} \le m_X \le \frac{2^{1-1/d} e (\Gamma(d+1))^{1/d}}{d}.
$$
For $d=4$, this gives $1.8554 < m_X < 2.5297$.

In $\M^d$, the proper time between two timelike points $\bf a$ and $\bf b$ is proportional to $({\rm vol}[{\bf a},{\bf b}])^{1/d}$, so
Theorem~\ref{thm:BB-scale} shows that again, up to a constant factor, the height of the interval provides a good approximation to the proper time
between the two points.  As explained in Bachmat~\cite{Bachmat}, this result extends to Poisson processes in Lorentzian manifolds with curvature: if we
take a Poisson process in an interval in such a manifold, then the longest chain will follow the maximal curve between the endpoints, and its length,
suitably scaled, will be a good approximation to the proper time.

It would be interesting to have better bounds on the constant $m_X$ when $X$ is an interval in $\M^d$: as far as we know, the bounds above of
Bollob\'as and Brightwell are the best known.  It would also be interesting to have some tight bounds, upper and lower, on the variance of $H_X(n)$.
A straightforward application of Talagrand's inequality gives that the standard deviation is at most of order $n^{1/2d}$, but the example of
the Baik-Deift-Johansson Theorem suggests that the actual deviation may be much smaller.  This issue is discussed by Bachmat in~\cite{Bachmat},
where it is suggested on the basis of heuristics and numerical evidence that the standard deviation in $\M^4$ might be of order~$n^{1/24}$.

The implications of Theorem~\ref{thm:BB-scale} for causal set theory were spelled out by Brightwell and Gregory~\cite{BGreg}, where they explained
how the height of $P[{\bf a},{\bf b}]$ could operate as a definition of the distance between two timelike points.  What about the distance between two
spacelike points ${\bf a}$ and ${\bf b}$?  In the underlying manifold, this distance is well-defined, and we might hope that some simple
parameter of the partial order approximates spacelike distance up to a constant factor.  For a causal set obtained from a Poisson process in $\M^2$,
one can define a concept of two points being ``nearest neighbours'', and then derive a distance function.  Brightwell and Gregory~\cite{BGreg} pointed out
a phenomenon that makes any such definition problematic in higher dimensions.  However exactly the definition is framed, we would
want that a {\em sufficient} condition for a spacelike pair ${\bf a},{\bf b}$ to be nearest neighbours is that there are points ${\bf c} < {\bf d}$
of the causal set such that ${\bf a}$ and ${\bf b}$ are the only two points of the causal set in $[{\bf c},{\bf d}]$.
However, given any spacelike pair ${\bf a}$, ${\bf b}$ in $\M^d$ ($d \ge 3$), by applying larger and larger boosts in a direction perpendicular to the line between ${\bf a}$ and ${\bf b}$, one can find an infinite collection of intervals in $\M^d$ containing ${\bf a}$ and
${\bf b}$ in the interior, all of essentially the same volume, each pair of intervals disjoint apart from an arbitrarily small tube around the line between
${\bf a}$ and ${\bf b}$.  Given such a collection, where the tube is small enough to contain no points of the causal set, almost surely one of these
intervals will contain exactly two points ${\bf c}$ and ${\bf d}$ of the causal set, with ${\bf a}$ and ${\bf b}$ in $[{\bf c}, {\bf d}]$.  Thus every
spacelike pair of points would be nearest neighbours.

The existence of these
rare ``voids'', where there are very few points of the Poisson process in an interval of large volume, along with the Lorentz invariance of the
manifold, dooms any attempt to define a distance in terms of minimising some quantity over some family of intervals.
It is, however, worth noting that this precise objection is of a theoretical nature only: Dowker, Henson and Sorkin~\cite{DHS} calculate that the probability
of a void of the size of an atom contained within the observed universe is negligibly small.

There remains no strong candidate for a natural 
definition of spacelike distance within a causal set.  See however Rideout and Wallden~\cite{RW} for some
practical approaches to measuring spacelike distances that yield useful answers.

Another direction of research concerns the characterisation of which finite partial orders can be embedded in $\M^d$ (i.e, of $d$-sphere orders), for
each~$d$.  For $d=2$, the answer is the family of 2-dimensional orders, and this is a well-studied and well-understood family.  For higher dimensions,
rather little is known.  Brightwell and Winkler~\cite{BW-sphere}, and independently Meyer~\cite{Meyer}, constructed, for each $d$, an example of a poset
that is a $d$-sphere order but is not a
$(d-1)$-sphere order.  Alon and Scheinerman~\cite{AS} showed that the number of $d$-sphere orders on $m$ labelled elements grows as $m^{md (1+o(1))}$ for
each fixed~$d$ (which provides a non-constructive proof of the result of Brightwell/Winkler and Meyer).  This number is comparable to the number of $m$-element
partial orders embeddable in $\R^d$, i.e., the number of $m$-element partial orders of dimension at most~$d$, and grows slowly compared to the total number
of $m$-element partial orders, which is $2^{m^2/4 + O(m)}$ -- we shall discuss this further later.  For a while, it was an open
question whether (a)~all 3-dimensional partial orders can be embedded in $\M^3$, and (b)~every partial order can be embedded in some~$\M^d$.  In 1999, Felsner,
Fishburn and Trotter~\cite{FFT} answered both questions by showing that there is some 3-dimensional partial order that cannot be embedded in any $\M^d$.
Their example is produced through Ramsey-theoretic arguments, and its size is enormously larger than the supposed number of points of spacetime in the
observable universe, but it is very likely that there are examples with fewer than a hundred elements: in~\cite{BW-sphere}, Brightwell and Winkler
suggested that the partial order $2^6$ of subsets of a 6-element set is not a sphere order in any dimension.
It is remarkable that this concrete question is apparently not easy to answer computationally; to the best of our knowledge, the complexity of determining
whether a given finite partial order is (for instance) a circle order is unknown.

If we are presented with a causal set, and told that it arises as a sample from a Poisson process on some Lorentzian manifold (perhaps with a few local
``defects''), can we determine the manifold, up to some rough isometry?  This is an ill-specified question, but it is nevertheless important for causal
set theory to have some kind of answer.

Bombelli and Meyer~\cite{BombM} showed that, if we are given a countable partial order, consisting of an iid sequence of points selected
according to the volume measure on some portion of a Lorentzian manifold of finite volume, then the manifold can be recovered from the
causal structure on the points.

If we are given instead a sequence of finite partial orders, sampled separately from a region of finite volume in a fixed manifold according to
Poisson processes with densities tending to infinity, then we would hope to be able to recover the manifold via a {\em poset limit}.
The general theory of poset limits has been developed recently by Janson~\cite{Janson} and
Hladk\'y, M\'ath\'e, Patel and Pikhurko~\cite{HMPP}; we give a brief account here.

Given finite posets $P$ and $Q$, the {\em density} $t(Q;P)$ of $Q$ in $P$ is the proportion of
$|Q|$-tuples of elements of $P$ that are isomorphic to $Q$.  A sequence $(P_n)$ of finite posets with $|P_n| \to \infty$ is said to be {\em convergent}
if, for every fixed finite poset $Q$, $t(Q;P_n)$ tends to a limit.  An appropriate limit object for a sequence of posets is what
has been informally dubbed a {\em poson}, by analogy with the idea of a graphon in the theory of graph limits.
Let $([0,1],\cB,\mu)$ denote the standard Borel probability space on the interval $[0,1]$, and $\prec$ the standard linear order on $[0,1]$.
A {\em poson} $W$ is a $\cB\times \cB$ measurable function from $[0,1]\times [0,1]$ to $[0,1]$ such that (i)~if $W(x,y) > 0$, then $x \prec y$, and
(ii)~if $W(x,y) >0$ and $W(y,z) > 0$, then $W(x,z) = 1$.  To generate a random $n$-element sample from a poson $W$, we take $n$ elements $X_1, \dots, X_n$
uniformly and independently from~$[0,1]$, and we put $X_i< X_j$ with probability $W(X_i,X_j)$, each pair independently.  Condition~(ii) then guarantees
that $<$ is a transitive relation, and condition~(i) ensures that the linear order induced by $\prec$ on the sample $(X_1, \dots, X_n)$ gives a linear
extension of the derived poset.  The density $t(Q;W)$ of a fixed poset $Q$ in a poson $W$ is the probability that a $|Q|$-element sample from $W$ is isomorphic
to $Q$; the density $t(Q;W)$ can be expressed as a multidimensional integral.  The sequence $(P_n)$ of posets is then said to converge to the poson $W$ if
$t(Q;P_n) \to t(Q;W)$ for every finite poset~$Q$.  Extending a result of Janson~\cite{Janson}, Hladk\'y, M\'ath\'e, Patel and Pikhurko~\cite{HMPP} proved
that every convergent sequence of posets converges to some poson.

The implications for causal set theory are indirect.  The theory of poset limits highlights the following conjecture of Bombelli~\cite{Bombelli}
as fundamental: if $U$ and $V$ are regions of Lorentzian manifolds of finite volume, such that there is no measure-preserving diffeomorphism
between the two, then there is some finite poset $Q$ such that $t(Q;U) \not= t(Q;V)$.  Here the density $t(Q;U)$ is the probability that an iid uniform
sample of $|Q|$ elements from $U$ yields a poset isomorphic to $Q$.

In causal set theory, we do not have a sequence of partial orders with sizes tending to infinity, but a single partial order that in most
ordinary terms would be thought of as very large (although not large enough to meaningfully apply the poset regularity lemma of~\cite{HMPP}), and the sense
in which it approximates a manifold will necessarily be quantitative, rather than a statement about limits.  One would like a result saying that, if
two finite portions of manifolds are both approximated well by the same partial order, then the two spaces are roughly isometric.  No specific conjecture
has been formulated, but this has nevertheless become known as the {\em Hauptvermutung} of causal set theory.

\section{Growth processes}

In the previous section, we discussed the properties of a Poisson process in a Lorentzian manifold.  However well such a model might fit observations,
it cannot make for a theory with any ``explanatory power'': we start with a manifold, embed the discrete set of points, and then ``rub out'' the unwanted
manifold at the end.  What would be far more satisfying is a random process, following some reasonably natural rules, whose samples typically resemble a
Poisson process in $\M^4$.

We know of no contender for such a process, and indeed it is surely too much to expect; any eventual theory will presumably not involve classical
probability measures on causal sets at all.  In the absence of a fully realised theory, what one can do is start from the other extreme, and ask what
form a suitable process might take, and then investigate processes of such a form.  If there are processes that share some of
the observed properties of spacetime, then all the better: if not, then one gains experience in asking the right questions and identifying the key
issues.

One outline structure that has arisen~\cite{RS-csg} is that of a {\em growth process}.  In this construction, points (of spacetime) are generated one at
a time, each new point being maximal among those seen so far, and placed above some subset of the existing points chosen according to some (random) rule.
In this way, an infinite causal set is generated.  The order in which the elements of the causal set is generated is part of the model, and it is
a {\em linear extension} of the causal set: a linear order on the ground-set extending the partial order.

One example of this process is known as the model of {\em random graph orders} in the mathematics literature, and {\em transitive percolation} in the
physics literature.  Here, we fix some value $p \in (0,1)$, and let the ground-set be $\{ x_1, x_2, \dots\}$.  For each $n$ in turn, a random
subset $S_n$ of $\{ x_1, \dots, x_n\}$ is chosen, with each element of $\{ x_1, \dots, x_n\}$ appearing in $S_n$ with probability~$p$, independent of all
other choices.  Then the new element $x_{n+1}$ is placed above all elements of $S_n$.  Typically, this rule will not obey transitivity, so additionally
$x_{n+1}$ is placed above all elements that are less than or equal to some element of $S_n$.  An alternative description of this random construction is
as follows: (i)~take a random graph on the ground-set, with each edge present independently with probability~$p$, (ii)~for $i<j$, interpret an edge from 
$x_i$ to $x_j$ as putting $x_i < x_j$, and (iii)~take the transitive closure.  
We may stop the process after $n$ elements have been added, in which case we denote the resulting random graph order $P_{n,p}$.
This model of random graph orders was first considered as a model of random acyclic directed graphs by Barak and Erd\H os~\cite{BE}, introduced in
its own right by Albert and Frieze~\cite{AF}, and then studied by a number of other authors since, including a series of papers by Bollob\'as and Brightwell~\cite{BB-rgo1,BB-rgo2,BB-rgo3}.

Obviously there is an enormous variety of other growth models, and it is difficult to say anything of interest about growth models in full generality.
Rideout and Sorkin~\cite{RS-csg} identified two extra conditions that they felt were desirable in a growth model representing the universe.
The first property they called ``general covariance'': for every $n$, and every poset $Q$ on $n$ elements, conditional
on the first $n$ points generated forming $Q$, each linear extension of $Q$ is equally likely to be the order in which the first $n$ points were generated.
One needs to be careful with this definition: the partial orders in Rideout and Sorkin's definitions are unlabelled, and they discuss their process in
terms of transition probabilities for a Markov chain whose states are unlabelled finite partial orders.  Given any unlabelled $m$-element partial order $Q$,
with $e(Q)$ linear extensions, each linear extension induces a path from the empty partial order to $Q$, adding one maximal element at a time.
Some paths may be induced by more than one linear extension, as there may be several isomorphic copies of some partial order that are stems of $Q$.
The condition of general covariance states that, for every $Q$, each path from the empty partial order to~$Q$ occurs with probability proportional to the
number of linear extensions of $Q$ inducing that path.

Rideout and Sorkin~\cite{RS-csg} dubbed their second property {\em Bell causality}.  It reflects the principle that a change in one part of the
``universe'' cannot affect some other part of the universe that is not in the future light cone of the change.  The translation of this to a growth model
is as follows.
A transition of the Markov chain involves adding a new maximal element: suppose we are considering two possible transitions from a state $Q$,
adding the new element above either of the two down-sets $S_1$ or $S_2$ of the existing elements, with transition probabilities $p(Q;S_1)$ and $p(Q;S_2)$ respectively.  Now consider transitions from a state $Q'$ which is the poset $Q$ restricted to $S_1 \cup S_2$: the condition of Bell causality is that
$$
p(Q;S_1) p(Q';S_2) = p(Q;S_2) p(Q';S_1).
$$
Thus the relative probabilities of the two transitions are unchanged on the removal of ``spectators'' which would be incomparable with the new element
in both transitions.

We observe that transitive percolation, for any value of the parameter~$p$, satisfies both of these conditions.  To see general covariance, it is enough to
show that, for any partial order $Q$ on $\{ x_1, x_2, \dots, x_m\}$ with $x_1 < x_2 < \cdots < x_m$ as a linear extension,
the probability that the first $m$ steps of the transitive percolation process produces $Q$ does not depend on the labelling of the elements of $Q$.
We see that the probability of $Q$ is given by $p^{c(Q)}(1-p)^{i(Q)}$, where $c(Q)$ is the number of covering pairs (i.e., pairs $(x,y)$ such that $x<y$
in $Q$ but there is no $z$ with $x < z < y$ in $Q$), and $i(Q)$ is the number of incomparable pairs.  Checking that transitive percolation satisfies Bell
causality is straightforward.

Rideout and Sorkin~\cite{RS-csg} also imposed an extra condition that, from every state, the probability of the new element being incomparable to every
existing element should be positive; we say that the growth model is {\em gregarious}.  They then proved that the only growth models satisfying general covariance, Bell causality and this extra condition are what they called {\em classical sequential growth models} or {\em csg models}.

A csg model is specified by a sequence of non-negative constants $t_0, t_1, \dots$, where $t_0 > 0$.  At any state $Q$, a partial order with $m$ elements,
a subset $S_m$ of the existing elements is selected, and the new element $x_{m+1}$ is placed above all elements below or equal to some element of $S_m$.
For each $S$, the probability that $S_m=S$ is proportional to $t_{|S|}$, so the probability that $S_m = S$
is $t_{|S|} \Big( \sum_{k=0}^m \binom{m}{k} t_k \Big)^{-1}$.  For instance, transitive percolation is the example where $t_k = (p/(1-p))^k$ for each~$k$.

To see that a csg model, as above, satisfies general covariance, we calculate the probability of generating a particular labelled partial
order on the set $\{x_1, \dots, x_n\}$, where $x_{i+1}$ is above a set $B_i$ of $b_i$ elements, and covers a set $C_i\subseteq B_i$ of $c_i$ elements.
Thus the given partial order is generated if each $S_i$ contains $C_i$ and is contained in $B_i$.  Now the probability that the given partial order is
generated is the product of ratios $a_i/\sum_{k=0}^i \binom{i}{k} t_k$, where $a_i = \sum_{\ell = c_i}^{b_i} \binom{b_i-c_i}{\ell-c_i} t_\ell$, the
total ``weight'' attached to all values of $S_i$ with $C_i \subseteq S_i \subseteq B_i$.
In the overall product, the denominators are purely functions of the $t_k$, while the numerator $a_i$ depends only on the structure of the partial order
below $x_{i+1}$.  Hence, on taking a different order of generation, the same denominators appear, and the same numerators appear in a different order,
so the ratio is unchanged.

The main theorem of Rideout and Sorkin~\cite{RS-csg} is as follows.

\begin{theorem} \label{thm:RS}
The only growth models that satisfy general covariance and Bell causality, and are gregarious, are csg models.
\end{theorem}

What if one drops the extra ``technical'' condition that the growth model be gregarious?  This allows also models of the following form: there are
a number of ``phases'' of the growth process, with each phase of the process developing
as a csg model and, with the possible exception of a final infinite phase, stopping after a pre-determined finite number of steps.
All elements of each new phase are above all elements of previous
phases.  The number of steps of the current phase, and the parameters of the csg model in effect during the phase, are arbitrary functions of the partial
order at the end of the previous phase.  Varadarajan and Rideout~\cite{VR} and Dowker and Surya~\cite{DS} showed that the growth models satisfying general covariance and Bell causality are exactly the processes of this form.  Thus dropping the requirement that our growth models be gregarious does not
introduce any different possible ``local'' structures.

Let us then consider whether csg models include some that could form the basis of a satisfactory theory of quantum gravity.  First of all, we consider
whether random graph orders might be a contender.  The large-scale structure of a random graph order $P_{n,p}$ has been extensively investigated, and is
now quite well understood.  Alon, Bollob\'as, Brightwell and Janson~\cite{ABBJ} showed that, for a fixed $p$ and an infinite ground-set, there are
a.s.\ infinitely many {\em posts}, which are elements of the partial order comparable to all others.  Sharper results about the distribution of the sequence
of posts were given by Bollob\'as and Brightwell~\cite{BB-rgo3} and Kim and Pittel~\cite{KP}.  For small $p$, posts become very rare: the probability that a given element is a post is asymptotically $2\pi e^{\pi^2/6} p^{-1} e^{-\pi^2/3p}$.  Bollob\'as and Brightwell~\cite{BB-rgo3} explain
how the presence of many posts implies that parameters such as the height of a random graph order are normally distributed, with mean proportional to~$n$.

If we instead take a function $p=p(n)$ such that $n p /\log p^{-1} \to \infty$ and $p \log n \to 0$, then the random graph order $P_{n,p}$ exhibits a
very different structure.  Suppose first that $n = \alpha p^{-1} \log(p^{-1})$, for $\alpha < 1$ a fixed constant.  Then the number of elements above
$x_1$ in the random graph order is fairly well concentrated around $n^\alpha$.  This implies that, in a random graph order, the element $x_m$ is
typically incomparable to most elements in the set $\{ x_k : |k-m| < p^{-1} \log (p^{-1}) \}$.  On the other hand, if $n \ge \alpha p^{-1} \log(p^{-1})$
for some $\alpha > 1$, and $I(1)$ denotes the number of elements incomparable to $x_1$, then it can be read out of a result of Pittel and
Tungol~\cite{PT} that, whenever $M = O(p^{-1} \log\log (p^{-1}))$,
$$
\Pr( I(1) < p^{-1}\log(p^{-1}) + M ) = (1+ o(1)) \exp( - e^{-pM}).
$$
This implies that each element $x_m$ is almost surely comparable to all elements
outside the range $\{ x_k : |k-m| < p^{-1} \log(p^{-1}) + \omega(n) p^{-1} \}$, where $\omega(n)$ is any function tending to infinity.
There are similar results in the earlier papers of Simon, Crippa and Collenberg~\cite{SCC} and Bollob\'as and Brightwell~\cite{BB-rgo2}.
For this range of parameters, Bollob\'as and Brightwell~\cite{BB-rgo1} showed that the width of the partial order (i.e., the size of a largest
antichain) is of order $p^{-1}$ and, moreover, a.s.\ every element is in an antichain of size of order $p^{-1}$.  By contrast, if we fix $p$ and let
$n$ tend to infinity, then there is the occasional antichain of size much larger than $p^{-1}$.

This is certainly not a picture that suggests a 4-dimensional Lorentzian manifold, and indeed it is possible to point out various concrete ways in which it
differs: for a start, almost surely every finite partial order will appear as an interval in a random graph order, whereas this is not the case in a Poisson
process in a manifold.

However, a random graph order does share some properties of interest with a Poisson process in a Lorentzian manifold.
Rideout and Sorkin~\cite{RS-limits} gave experimental evidence that a suitably scaled sequence of finite random graph orders approaches a
poset limit.  Namely, they looked at the densities $t(Q,P_{n,p(n)})$ of various small posets $Q$ in a random graph order $P_{n,p(n)}$, and observed that,
if $p(n)$ is tuned so that the density of comparable pairs remains constant, then the densities of all the other posets $Q$ they looked at also
approached limits.  Moreover, the densities of some suborders appeared to tend to~0.  They pointed out that these properties were what would be
expected on taking a sequence of Poisson processes, with increasing intensity, in a manifold.

\begin{figure} [hbtp]
\epsfxsize290pt
$$\epsfbox{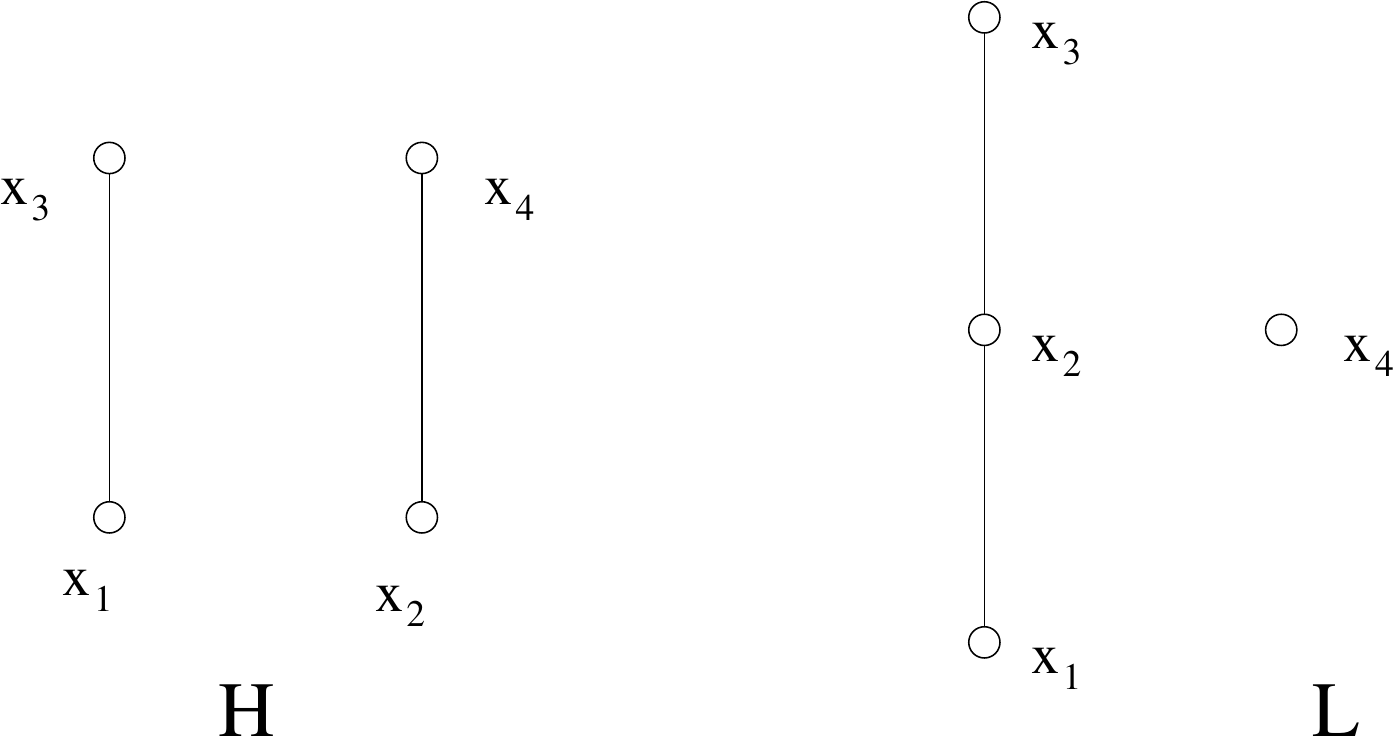}$$
\caption{The partial orders $H$ and $L$} \label{fig:HL}
\end{figure}

Brightwell and Georgiou~\cite{BGeorg} explained this phenomenon in terms of the known large-scale structure of a random graph order, as described above.
In the language of poset limits, as set out the previous section, they showed that a sequence $(P_{n,p(n)})$ of random graph orders
is convergent if and only if either $p^{-1} \log(p^{-1}) / n$ tends to a limit in $[0,1)$, or $\liminf p^{-1} \log(p^{-1}) /n \ge 1$.  In the case where $p^{-1}\log(p^{-1})/n$
tends to zero, the proportion of comparable pairs of elements in $P_{n,p(n)}$ tends to~1, as described above, and the poset limit is a chain.  In the
case where $\liminf p^{-1} \log(p^{-1}) /n \ge 1$, the random graph order is very sparse, and the limit is an antichain.  In the intermediate case, where
$p^{-1} \log(p^{-1}) /n \to c$, for $c \in (0,1)$, one can take the limit poson $W : [0,1]\times[0,1] \to [0,1]$ to be given by:
$$
W(x,y) = \begin{cases} 1 & \mbox{if } y > x+c \\ 0 & \mbox{if } y \le x+c \end{cases}.
$$
The finite partial orders with positive density in $W$ are all {\em semiorders}: they can be represented by a family of unit-length intervals in $\R$ with
$x<y$ if the interval representing $x$ lies entirely to the left of the interval representing~$y$.
An alternative characterisation is that semiorders are partial orders containing neither of the posets $H$ or $L$ in Figure~\ref{fig:HL} as suborders.
It is natural to extend the definition of semiorders to include posons; a poson $W$ is a semiorder if $W(x,y)$ is 0-1 valued, non-increasing in $x$ and
non-decreasing in~$y$.  Chains and antichains are semiorders, so the limit of any convergent sequence $(P_{n,p(n)})$ of random graph orders is a semiorder.

The causal structure of the observed universe is manifestly not that of a semiorder; a semiorder has a very well-defined temporal order, but no spatial
structure at all.  The conclusion is that a random graph order, even with a carefully chosen parameter, is not an adequate model.

Might it be that general csg models exhibit very different behaviour?  There are certainly some interesting cases where the random partial order
produced is sparse, for instance if all the parameters $t_k$ are zero beyond a certain point.  The case where $t_0 , t_1 > 0$ and $t_2=t_3=\cdots =0$ produces
a random forest: each new element $x_m$ is placed above at most one other, with $x_m$ being minimal with probability $t_0/(t_0 + (m-1)t_1)$, and placed
above each existing element with probability $t_1/(t_0 + (m-1)t_1)$.  In the limit, there are a.s.\ infinitely many connected components, and every element
is covered by infinitely many later elements.  The next case is that of ``random binary orders'', where each arriving element is placed above (up to) two existing elements, chosen uniformly at random.  The structure of random binary orders was investigated by Georgiou~\cite{Georg}.  One of the main results of
that paper is that, for any constants $\eps, \eta \in (0,1)$, there is a constant $M$ such that, for each fixed $m$, with probability at least $1-\eta$,
the number of elements incomparable with $x_m$ is at most $M m^{2+\eps}$, and $x_m$ is below all elements $x_n$ with $n > m^{4+\eps}$.  There is some
delicacy in the statement, as there will be a few exceptional elements that are comparable with very few elements for an unusually long time in the
growth process, and indeed there are a.s.\ infinitely many $m$ such that $x_m$ is incomparable to at least about $m^3$ elements.  Georgiou~\cite{Georg}
proves that, a.s., for every $m$, $x_m$ is incomparable with at most $m^{27/5}$ elements, and is below every $x_n$ with $n \ge m^{12}$.  The basic
phenomenon here is that, even in this csg model, which is as sparse as possible without being the model of random forests, every element is incomparable
with finitely many others.  One might conceivably be interested in this behaviour in the context of a model of a very rapidly expanding universe.

However, for our purposes it seems more relevant to consider ``dense'' csg models, where the proportion of comparable pairs is positive.  For
comparison purposes, we should consider a sequence $(C_n)$ of finite partial orders, with $C_n$ equal to the induced csg model on $\{ x_1, \dots, x_n\}$
with parameter sequence ${\bf t}(n)$.  Brightwell and Georgiou~\cite{BGeorg} showed that, even with this extra generality, with an infinite collection
of parameters instead of a single one, the only possible poset limits for a sequence of csg models are semiorders.

\begin{theorem} \label{thm:semiorder}
Let $(C_n)$ be a sequence of finite csg models, where $C_n$ has $n$ elements for each~$n$.  If the sequence $(C_n)$ is convergent,
then the limit is a semiorder.
\end{theorem}

To prove Theorem~\ref{thm:semiorder}, Brightwell and Georgiou showed directly that, for any sequence of finite csg models with sizes tending to
infinity, the expected numbers of copies of the posets $H$ and $L$ in Figure~\ref{fig:HL} tend to~0.

Moreover, Brightwell and Georgiou~\cite{BGeorg} determined conditions under which a sequence of csg models does converge to a poset limit.  The basic idea
is that, in any csg model, the structure is locally quite similar to that of a random graph order.  The key sequence of parameters is $(p_m)$,
where $p_m = \E |S_m| /m$, the probability that any given element $x_n$, with $n < m$, is selected for the set $S_m$; at the introduction of the new element
$x_{m+1}$, $p_m$ plays a role analogous to the parameter~$p$ of random graph orders.  The sequence $(p_m)$ of probabilities can take arbitrary
jumps upwards within $(0,1)$: the parameter $t_m$ of the csg model plays no role in generating the sets $S_1, \dots, S_{m-1}$,
but, if it is large enough, then there will be a high probability that $S_m = \{ x_1, \dots, x_m\}$, and so $x_{m+1}$ will be placed above all
existing elements.
Brightwell and Georgiou showed that $\E |S_m|$ is non-decreasing in $m$, so the sequence $(p_m)$ can only decrease slowly.  Therefore, in all
interesting cases, $(p_m)$ is slowly varying over long stretches and, if $p_m$ is close to a value $p$ over a long stretch, then the structure of the csg
model is, on a large scale, very close to that of a random graph order with parameter~$p$.

Negative results such as Theorem~\ref{thm:semiorder} were expected, and there remains interest among physicists in the behaviour of random graph orders,
and more generally in csg models.  If there is to be a ``successful'' theory of causal sets, it will not be based on a csg model, indeed most likely not
on any classical stochastic process, but nevertheless a csg model may serve as a useful
toy model, having some features of the real theory without all of the complexity that will surely prove necessary.
For instance, one question that arose in the study of csg models was the nature of the appropriate $\sigma$-field of ``observable'' events.  The description
of the growth process entails a natural $\sigma$-field $\cG$, and an event $G$ in $\cG$
is deemed to be ``observable'' if it is label-invariant, i.e., if renumbering the elements of the partial order with a different linear extension does not
affect membership in~$G$.  Brightwell, Dowker, Garc\'ia, Henson and Sorkin~\cite{BDGHS} showed that, up to sets of measure zero, the $\sigma$-field of
observable events is generated by the {\em stem} events, which are the events of the form ``the causal set contains a down-set isomorphic to $L$'', for
each fixed finite partial order~$L$.

One feature that does have resonance for theories of quantum gravity is the existence of posts: recall that, for any fixed $p$, there are a.s.\
infinitely many posts in the random graph order with parameter~$p$.  This might represent the
possibility of a succession of universes, with no information flowing from one to the next.  The question then arises of which csg models have posts: this
is currently an open problem, but it is possible to suggest an answer.  If the sequence $(t_k)$ of parameters of the model behaves as $c^k$, for any
constant~$c$, whether smaller or greater than~1, then the csg model behaves as a random graph order on large scales.  If the sequence tends to infinity more
rapidly than an exponential function, then on large enough scales the csg model is denser than any random graph order, and the proportion of the elements
of $\{ x_1, \dots, x_n\}$ that are posts will tend to~1 as $n$ tends to infinity.  It is more interesting to look at sequences of parameters $(t_k)$
that tend to zero more rapidly than any exponential function.  For well-behaved sequences, the results of Brightwell and Georgiou~\cite{BGeorg} imply that
the csg model
then resembles a random graph order on any scale, but that as the scale increases the effective parameter $p_n$ of the random graph order decreases,
implying that posts become rarer as the number of elements increases.  Calculations strongly suggest that the threshold for the existence of
infinitely many posts is $p_m \simeq \pi^2 / 3 \log m$.  Thus we expect that, if $t_k \simeq (t/\log k)^k$ for some $t < \pi^2/3$, then there are
a.s.\ finitely many posts, while, if $t_k \simeq (t/\log k)^k$ for some $t > \pi^2/3$, then there are a.s.\ infinitely many posts.  Establishing a result
along these lines is work in progress.

Sorkin's idea has always been that the eventual theory of causal sets will not involve a model using ``classical'' probability measures at all.
Quantum mechanics can be described using quantum measures, obeying some but not all of the rules of probability measures, and it is thus at least
worth considering the possibility that the evolution of spacetime also involves quantum measures.

Another possibility before giving up on classical probability altogether is to challenge one of the hypotheses of Theorem~\ref{thm:RS}, either
general covariance or Bell causality.  Of the two, general covariance seems
mathematically more natural, and also to have more unassailable motivation: it embodies the principle that if we observe the spacetime universe,
we see no effects that can be attributed to a specific favoured linear extension of the partial order.  For Bell causality, while the principle is
appropriate for events taking place in the universe, it does not seem automatic that it should apply to the generation of the universe itself.
In the next section, we explore what happens if we drop the condition of Bell causality.

\section{Order-invariance}

In this section, we give a short account of the work in two papers~\cite{BL1,BL2} of the authors.  The idea behind these papers is to examine the property
of general covariance, introduced in the previous section, in its own right.  As we shall explain, measures satisfying general covariance are
analogues of the Gibbs measures of statistical physics, and we believe that the theory we develop here will prove to be of interest beyond any
relevance it has to the causal set programme.

As mentioned earlier, Rideout and Sorkin~\cite{RS-csg} expressed the property of general covariance in terms of transitions of a Markov chain whose
states are {\em unlabelled} finite partial orders.  It is very often cleaner to work with labelled combinatorial structures, and accordingly we shall
discuss Markov chains whose states are labelled partial orders, and in this context we use the term ``order-invariance'' instead of ``general covariance'':
we give a definition later.  As with growth models, we start with the empty poset, and at each stage we add one
new maximal element, keeping track of the order in which the elements are generated.
The infinite poset $P=(Z,<)$ generated is always a causal set, and the order in which the elements are generated is a {\em natural extension} of $P$,
i.e., a bijection from $\N$ to $Z$ whose inverse is order-preserving.  Alternatively, we can view a natural extension of $(Z,<)$ as an enumeration of
$Z$ that respects the partial order~$<$.

We begin with an example of a very different nature from that of csg models.  Here, we start with a fixed causal set $P=(Z,<)$, and give a prescription
for generating a random natural extension of $P$.  Our process will have the property of order-invariance: if we stop after some fixed number of steps, then,
conditioned on the structure of the causal set (i.e., the elements included and the comparabilities among them), every possible order of generation is
equally likely.

\begin{figure} [hbtp]
\epsfxsize290pt
\vspace{-3.3truein}
$$\epsfbox{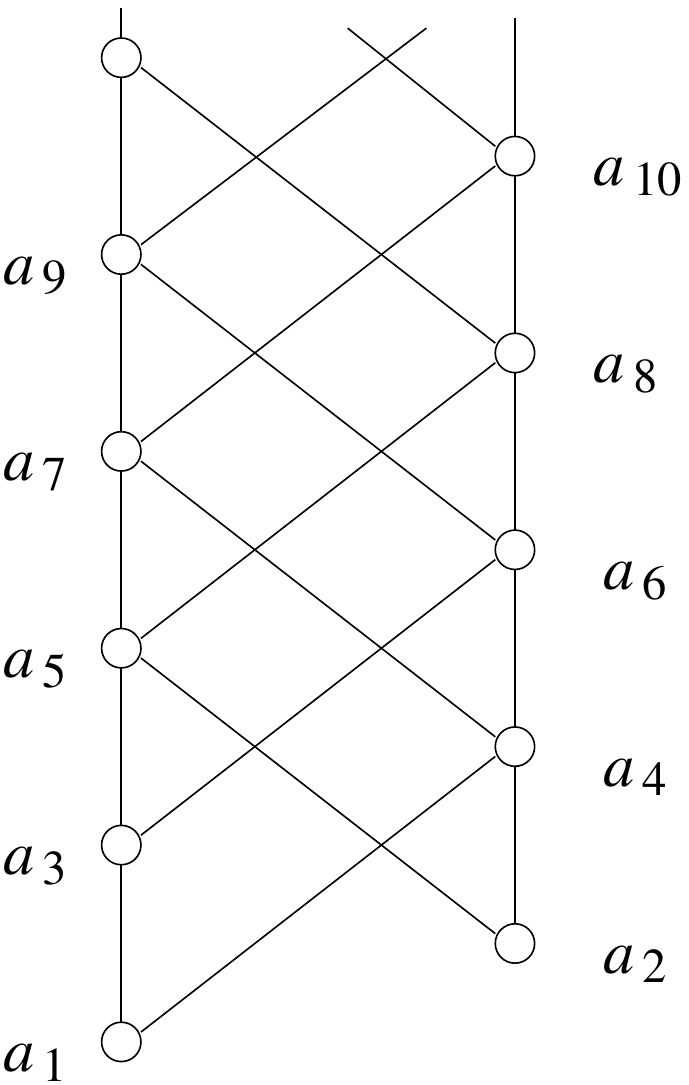}$$
\caption{The causal set $P=(Z,<)$} \label{fig:ladder}
\end{figure}

Consider the poset $P=(Z,<)$ in Figure~\ref{fig:ladder}.

Here $Z=\{a_1,a_2, \ldots \}$, and $a_j > a_i$ if $j > i+1$.
A bijection $\lambda: \N \to Z$ is then a natural extension if $X_k=\{\lambda (1), \ldots, \lambda (k)\}$ is a down-set in $P$ for each~$k$.
We seek a probability measure $\mu$ on the set $L(P)$ of natural extensions of $P$ which is ``uniform'' (in a suitable sense).

For this example, the $\sigma$-field of events is generated by events of the form
$$
E(a_{i_1}a_{i_2} \cdots a_{i_k}) = \{\lambda: \lambda (j) = a_{i_j} \mbox{ for } j=1, \ldots, k\},
$$
the set of natural extensions with ``initial segment'' $a_{i_1}a_{i_2} \cdots a_{i_k}$, for $k \in \N$ and the $i_j$ distinct elements of $\N$.
We call $a_{i_1}a_{i_2} \cdots a_{i_k}$ an {\em ordered stem} if $\{a_{i_1},a_{i_2}, \ldots, a_{i_j}\}$ is a down-set in $P$, for $j=1, \ldots, k$: i.e.,
if there is a natural extension of $P$ with this initial segment.  Our order-invariance condition is that, whenever
$a_{i_1}a_{i_2} \cdots a_{i_k}$ and $a_{\ell_1}a_{\ell_2} \cdots a_{\ell_k}$ are ordered stems such that
$\{a_{i_1}, a_{i_2}, \ldots, a_{i_k}\} = \{a_{\ell_1}, a_{\ell_2}, \ldots, a_{\ell_k}\}$, we have
\begin{equation} \label{eq:1}
\mu (E (a_{i_1}a_{i_2} \cdots a_{i_k})) = \mu (E (a_{\ell_1}a_{\ell_2} \cdots a_{\ell_k})).
\end{equation}

We describe the measure $\mu$ via a random process for generating $\lambda (1), \lambda (2), \ldots $ sequentially.
Given the set $X_k=\{\lambda (1), \ldots, \lambda (k)\}$, the element $\lambda (k+1)$ has to be one of the minimal elements of $P \setminus X_k$, and
there are at most two of these.  Our process is defined by the following rules:
if there is only one minimal element $a_k$ of $P \setminus X_k$, then we take $\lambda (k+1) = a_k$ with probability~1;
if there are two minimal elements $a_{k+1}$ and $a_{k+2}$ of $P \setminus X_k$, then we set $\lambda (k+1) = a_{k+1}$ with probability $\phi$ and
$\lambda (k+1) = a_{k+2}$ with probability $1-\phi$, where $\phi^2 = 1-\phi$.

It is easily seen that $\mu (E(a_2a_1)) = 1-\phi$: we choose $\lambda (1) =a_2$ with probability $1-\phi$, and then $a_1$ is the only
minimal element of $P \setminus \{a_2\}$, so we choose $\lambda (2) = a_1$ with probability $1$.
Also, $\mu (E(a_1a_2)) = \phi^2$, since we choose $\lambda (1) =a_1$ with probability $\phi$, and then $\lambda (2) =a_2$ with probability $\phi$.
Hence the condition of order-invariance forces the choice of $\phi$ to satisfy $\phi^2 = 1-\phi$.
One can extend the argument above to show that the random process defined above satisfies all the identities of the form~(\ref{eq:1}), so that
the process is order-invariant.

We have described a process generating a random natural extension of a fixed $P$, but we can also think of it as a growth process, growing a causal
set by adding one new maximal element at each step: we always get the same infinite causal set $P$, but the order in which the elements are generated
is random.

There is an alternative way to obtain the measure $\mu$ in our example.  Consider the sets $Z_n = \{a_1, \dots, a_n\}$,
the finite posets $P_n = P_{Z_n}$ induced by $P$ on each $Z_n$, and the uniform measures $\nu^{P_n}$ on
their sets of linear extensions, for each $n$.  It can be shown that
$$
\nu^{P_n}(E(a_{i_1}a_{i_2} \cdots a_{i_k})) \to \mu(E(a_{i_1}a_{i_2} \cdots a_{i_k})),
$$
as $n \to \infty$, for each ordered initial segment $a_{i_1}a_{i_2} \cdots a_{i_k}$.

For another example, we consider growth processes where the causal set generated is a.s.\ an
antichain.  If we require our causal sets to be labelled, then a growth process which
a.s.\ generates an antichain is nothing more than a sequence of random
variables: the labels of the elements, in the order they are introduced.
Order-invariance in this case requires that, if we condition on the set of the first
$k$ labels, for any $k$, then each of the $k!$ orderings of these labels
is equally likely.  That is, the sequence of labels must be {\em exchangeable}.

One way to generate a sequence of exchangeable random labels is to take
any probability distribution $\tau$ on any set $X$ of potential labels,
and let the labels be an iid sequence of random elements of $X$ with
probability measure~$\tau$.  Our labels are to be a.s.\ distinct,
so the probability measure $\tau$ must be atomless.
By the Hewitt-Savage Theorem~\cite{HS}, every sequence of
exchangeable random variables is a mixture of sequences of the type
above (i.e., there is a probability measure $\rho$ on some space
of probability measures on a set $X$: one measure $\tau$ is chosen
according to $\rho$, and then an iid sequence of random elements of $X$ is
generated according to $\tau$).

A similar mechanism allows us to incorporate csg models into our labelled context.  At each step, as a new element arrives, we
generate the down-set below it as in the previous section, and we equip it with a label chosen uniformly, and independently of everything else,
from $X=[0,1]$, with its usual Borel $\sigma$-field and Lebesgue measure.
We make this, fairly arbitrary, choice of $X=[0,1]$ as a set of potential labels for all our processes.

We now give the general definition of order-invariance.
We specify the outcome of a growth process
by giving an infinite string of (labels of) elements,
$x_1x_2\cdots$ ($x_i \in [0,1]$) in the order of their generation, together with a suborder
$<^\N$ of the index set $\N$ with its standard order: $i<^\N j$ if and
only if $x_i < x_j$ in the causal set $P=(X,<)$ generated.  Let
$\Omega$ be the set of pairs $\omega = (x_1x_2 \cdots, <^\N)$, the outcome space.
Let $\pi_k$ be the projection map on $\Omega$ given by $\pi_k (x_1x_2 \cdots, <^\N) = (x_1 \cdots x_k, <^{\N}_{[k]})$, where $<^{\N}_{[k]}$ is the
restriction of the order $<^\N$ to $[k]$.
Then we can think of $\pi_k (\omega)$ as the state of the process at step $k$, corresponding to outcome $\omega$.

Our basic events are subsets of $\Omega$ of the form $E(B_1, \cdots, B_k, <^{[k]})= \pi_k^{-1} (B_1, \cdots, B_k, <^{[k]})$, where each
$B_i$ is a Borel set and $<^{[k]}$ is a suborder of the natural order on $[k]$. We let ${\mathcal F}_k$ be the $\sigma$-field generated by the sets
$E(B_1, \cdots, B_k, <^{[k]})$, and ${\mathcal F}$ be the $\sigma$-field generated by $\cup_{k=1}^{\infty} {\mathcal F}_k$.
A growth process then gives rise to a probability measure $\mu$ on $(\Omega, {\mathcal F})$, which we call a {\em causal set measure}.
A causal set measure $\mu$ is {\em order-invariant} if, for any Borel sets $B_1, \ldots, B_k$, any suborder $<^{[k]}$ of $[k]$, and any linear extension $\lambda$ of $<^{[k]}$, we have
$$
\mu (E(B_1, \cdots, B_k, <^{[k]})) = \mu (E(B_{\lambda (1)}, \cdots, B_{\lambda (k)}, \lambda [<^{[k]}])),
$$
where the linear order $\lambda [<^{[k]}]$ on $\{ 1, \dots, k\}$ is given by $i (\lambda [<^{[k]}]) j$ if and only if
$\lambda (i) <^{[k]} \lambda (j)$.  This definition incorporates the notion of general covariance from the previous section, as well as that of
an order-invariant measure on the space of natural extensions of a fixed order, as in our earlier example.

Ideally we would like to classify order-invariant measures, in a manner similar to Theorem~\ref{thm:RS} of Rideout and Sorkin.  In reality, this is
too ambitious a goal, but some partial steps in this direction are taken in~\cite{BL1}.  It is noted in that paper that the space of order-invariant
measures is a convex subset of the set of all probability measures on $(\Omega,{\mathcal F})$.  This means that an arbitrary mixture of order-invariant measures
is also order-invariant, and suggests a focus on the {\em extremal} order-invariant measures: those that cannot be written as a convex
combination of two others.

An order-invariant measure that almost surely produces one fixed causal
set $P=(Z,<)$, as in our first example, will be called an order-invariant measure
on $P$. The process in that example is the only order-invariant
measure on that poset $P$, and it is extremal. This is however not true in general: there are causal sets with no order-invariant measures on them,
and others with infinitely many extremal order-invariant measures.

An example of a causal set admitting no order-invariant measures is an infinite labelled antichain; such an order-invariant measure
would give every element of the antichain the same probability of appearing as the bottom element in a random natural extension, and
this is not possible.  On the other hand, from the Hewitt-Savage Theorem, extremal order-invariant measures a.s.\
giving rise to an antichain are effectively the same as iid sequences of random elements of $[0,1]$.

The main result of~\cite{BL1} shows that all extremal order-invariant measures on $(\Omega, {\mathcal F})$ are, in a sense, a combination of
extremal order-invariant measures of these two types.

\begin{theorem} \label{thm:extremal}
Let $\mu$ be an extremal order-invariant measure.  Then there is a poset
$Q=(Z,<)$, either a causal set or a finite poset, with a marked set $M$ of
maximal elements such that, if $Q'$ is obtained from $Q$ by replacing each
element $z$ of $M$ with a countably infinite antichain $A_z$, then the
poset $\Pi$ generated by $\mu$ is a.s.\ equal to $Q'$, except for the
labels on the antichains~$A_z$.
\end{theorem}

Probably the most interesting special case is when the set $M$ of marked
maximal elements is empty, so that the extremal measure $\mu$ is an
order-invariant measure on the fixed (labelled) causal set $Q$.
However, not every order-invariant measure on a fixed causal set is extremal, so the theorem falls short of characterising extremal order-invariant measures.
The other extreme case is when $Q$ consists of a single marked element
$z$, and so $Q'$ consists of the single antichain $A_z$
(corresponding to an exchangeable sequence of random labels $[0,1]$ from an atomless
probability distribution on $[0,1]$).

The concept of order-invariance on fixed causal sets has occurred elsewhere in the literature, often in different guises.
Brightwell~\cite{B88} studied random linear extensions of locally finite posets.  His main theorem
gives that, if, for some fixed $t$, a causal set $P$ is such that every element is incomparable
with at most $t$ others, then there is a unique order-invariant measure on~$P$.
The case when $P$ is the two-dimensional grid $(\N \times \N,<)$ is linked to
representation theory and harmonic functions on the Young lattice (the lattice of down-sets of the grid) -- see
Kerov~\cite{Kerov}, for instance.  The topic of order-invariant measures on fixed causal sets is considered in detail in~\cite{BL2},
where there are many more examples and results.

The family of natural extensions of a fixed causal set $P$ can also be
viewed as the set of configurations of a (1-dimensional) spin system, and
order-invariant measures can then be interpreted as Gibbs measures,
so that some of the general results discussed in, e.g., Georgii~\cite{Georgii} apply.
Specifically, we now give an alternative characterisation of order-invariance.

For $<^\N$ a suborder of $\N$, $k \in \N$, and $\lambda$ a linear
extension of $<^\N_{[k]}$, let $\lambda^+$ be given by
$$
\lambda^+(i) = \begin{cases} \lambda(i) & i \le k \\ i & i > k.
\end{cases}
$$
So $\lambda^+$ is the natural extension of $<^\N$ obtained from the natural order on $\N$
by permuting the first $k$ natural numbers according to $\lambda$.
Now, for $\omega = (x_1x_2\cdots, <^\N) \in \Omega$, $k \in \N$, and $\lambda$ a linear extension
of $<^\N_{[k]}$, set $\lambda^+[\omega] = (x_{\lambda^+(1)}x_{\lambda^+(2)}\cdots, \lambda^+[<^\N])$;
this is the element of $\Omega$, representing the same labelled poset, obtained by permuting the first~$k$ elements of the 
ordered stem according to~$\lambda$.  Now, for $E \in \mathcal F$, let $\nu^k(E)(\omega)$ be the proportion of linear
extensions $\lambda$ of $<^\N_{[k]}$ such that $\lambda^+[\omega] \in E$.

For any $\omega \in \Omega$ and any $k$, the function $\nu^k(\cdot)(\omega)$
gives a probability measure on $\mathcal F$, namely the uniform measure on
$\lambda^+[\omega]$, where $\lambda$ runs over linear extensions of~$<^\N_{[k]}$.
Now we have our alternative characterisation of order-invariance.

\begin{theorem} \label{thm:DLR}
Let $\mu$ be a causal set measure.  Then $\mu$ is order-invariant if and only
if
$$
\mu(E) = \E_\mu \nu^k(E)
$$
for every $E\in \mathcal F$ and every $k\in \N$.
\end{theorem}

This is an analogue of the DLR equations from statistical physics
about conditional probabilities.  It corresponds to specifying a
boundary condition outside a finite volume -- here this means
that we condition on all the information about $\omega$ except the order
in which the first $k$ elements are generated.
We then realise the conditional Gibbs measure, which in our setting is $\nu^k(\cdot)(\omega)$.

We now describe a sequence $(\cG_k)_{k \in \N}$ of $\sigma$-fields, which correspond to the {\em external $\sigma$-fields} in interacting particle systems.
For $k \in \N$, ${\mathcal G}_k$ is, roughly speaking, the family of sets in $\mathcal F$ that are
invariant under all permutations of the first $k$ elements producing another natural extension.  We also define the {\em tail $\sigma$-field}
${\mathcal G} = \bigcap_{k=1}^\infty {\mathcal G}_k$.

The following result, again from~\cite{BL1}, establishes that the functions $\nu^k (\cdot) (\cdot)$ are a family of measure
kernels, analogous to Gibbsian specifications in statistical mechanics.

\begin{theorem} \label{thm:nu-k}
Let $\mu$ be an order-invariant measure.  For any event $E \in \mathcal F$, and
any $k \in \N$, we have
$$
\mu(E \mid {\mathcal G}_k) = \nu^k(E),
$$
almost surely.
\end{theorem}

We are now in a position to describe some conditions on an order-invariant measure~$\mu$ that are equivalent to extremality.
In general, the limit of the measures $\nu^k$ is not necessarily equal to $\mu$, but we do have the following.

\begin{theorem} \label{thm:nu}
Let $\mu$ be any order-invariant measure, and let $E$ be any event in
$\mathcal F$.  Then the sequence $\nu^k(E)(\omega)$ converges to a
$\mathcal G$-measurable limit $\nu(E)(\omega)$, $\mu$-a.s.  Moreover,
$\nu(E) = \mu(E \mid {\mathcal G})$, $\mu$-a.s., and $\E_\mu \nu(E) = \mu(E)$.
\end{theorem}

We say that an order-invariant measure $\mu$ is {\em essential} if, for
every $E \in {\mathcal F}$, $\nu^k(E) \to \mu(E)$ a.s.  In other words, $\mu$ is
essential if, for every $E$, the limit $\nu(E)$ in Theorem~\ref{thm:nu} is
a.s.\ equal to $\mu(E)$, or equivalently $\mu(E \mid {\mathcal G}) = \mu(E)$ a.s.
We say that $\mu$ has {\em trivial tails} if $\mu(H)$ is equal to~0 or~1
for every $H \in {\mathcal G}$.
It is shown in~\cite{BL1} that, for any order-invariant measure $\mu$, the following are equivalent:
(i)~$\mu$ is extremal, (ii)~$\mu$ is essential, (iii)~$\mu$ has trivial tails.

We finish this section by returning to the application to causal set theory.
We saw in the previous section that there is no csg model whose samples resemble those of a Poisson process in $\M^4$, or in any other manifold.
Now that we have relaxed the conditions on our growth processes, can we find an order-invariant causal set process that produces a random causal set
whose distribution is that of a Poisson process in (to be specific) the future light cone of a point in $\M^4$?  This question was raised in~\cite{BL2},
and it remains an open problem.   We discuss below the presumably simpler question with $\M^2$ (or, equivalently, $\R_+^2$ with the
coordinate order) replacing~$\M^4$, but the same principles apply to $\M^4$.

Suppose we take a Poisson process $X$ in the positive
quadrant $\R^2_+$, and consider the induced poset $P=(X,<)$: if we can show that there is a.s.\ an invariant measure on $(X,<)$, then we
can represent the overall measure as a mixture of these measures.
We now take a sequence of larger and larger squares $S_N = [0,N]\times[0,N]$, and consider the finite partial orders $P_N = P_{X \cap S_N}$, the
restrictions of $P$ to the points of the Poisson process in the squares $S_N$.  A natural candidate for an order-invariant measure $\mu$ on~$(X,<)$
is the limit of the uniform measures $\nu^{P_N}$, as earlier.  Now we focus on one minimal element $x$, say the one closest to the origin in Euclidean
distance, and consider $q_N(x) = \nu^{P_N}(E(x))$: the probability, in a uniform random linear extension of $P_N$, that $x$ is the bottom element.
If the limit of measures exists, then the sequence $(q_N(x))$ of probabilities will have to tend to a non-zero limit $q(x)$; indeed, the sum of the
$q(x)$ over all minimal elements $x$ will have to be equal to~1.  On the other hand,
if $q(x) = \lim_{N \to \infty} q_N(x)$ exists and is non-zero, then presumably the same is true of other
similar limits, and we should be able to use the values of these limits to define an order-invariant measure on $(X,<)$ (to begin, at the first step
of the growth process, we would select each minimal element $x$ with probability $q(x)$).
So the question of whether an order-invariant measure on $(X,<)$ exists would seem to boil down to the question of whether the sequence $(q_N(x))$ a.s.\
tends to a non-zero limit.  Unfortunately, our ability to estimate numbers of linear extensions for the partial orders $P_N$ is far too limited to
enable a direct attack on this question.

\section{Uniform measures}

A ``path integral'' or ``sum over histories'' approach to quantum gravity would proceed by specifying a space of ``potential'' histories, and assigning an
amplitude (or quantum measure) to sets of those histories.  The idea is that most of the ``weight'' of the measure is on histories similar to the observed
universe.  In other words, similarity to a Poisson process in a 4-dimensional Lorentzian manifold would be
an ``emergent property'' from some simple and natural random process.  A csg model constitutes one example of such an approach, where the probability measure
is classical.  Another approach is to examine what happens when we take very simple probabilistic models, such as the uniform model on (some class of)
partial orders of fixed size.

If we fix the number of elements at some (large) value $n$, then a uniformly random partial order on $n$ elements is almost surely a 3-layer poset.
A {\em $k$-layer poset} is one where the elements are partitioned into layers $A_1, A_2, \dots, A_k$, such that (i)~whenever $x\in A_i$ is below $y \in A_j$
in the poset, we have $i<j$, and (ii)~all elements of $A_i$ are below all elements of $A_{i+2}$, for each $i=1, \dots, k-2$.  The Kleitman-Rothschild
Theorem~\cite{KR} states that almost all $n$-element partial orders are 3-layer posets in which $|A_1|$ and $|A_3|$ are about $n/4$, and $|A_2|$ is about $n/2$.
Notice that, given a partition into sets of the appropriate sizes, and any subset $S$ of
$(A_1 \times A_2) \cup (A_2 \times A_3)$, we may form a 3-layer poset by putting $x<y$ only for those pairs in $S \cup (A_1 \times A_3)$.  Thus the number
of 3-layer partial orders is approximately $2^{n^2/4}$.

Brightwell, Pr\"omel and Steger~\cite{BPS} gave the following sharper estimate for the number $C_n$ of partial orders on $n$ labelled elements: for some
absolute constant $C>1$,
$$
C_n = (1+O(C^{-n})) \sum_{s=0}^{n} \binom{n}{s} 2^{(s+1)(n-s)} = (1+O(1/n)) \beta 2^{(n+1)^2/4} \binom{n}{\lfloor n/2 \rfloor}.
$$
In the sum, the $s$-term is a very tight estimate of the number of 3-layer posets with $s$ elements in the middle layer, and the result implies that
the proportion of $n$-element posets that are not 3-layer posets is exponentially small.
In the second expression, the term $\beta$ takes one constant value if $n$ is even and a very slightly different one if $n$ is odd.

A typical 3-layer poset has no ``spatial structure'', and does not resemble a Poisson process in any manifold, so it is unlikely that the uniform model has
much direct relevance to the causal set programme.  What might however be of interest is to restrict attention to certain classes of partial orders, and/or to
assign different ``activities'' to different $n$-element partial orders, as is common in statistical physics models.

One simple step is to constrain the number of comparable pairs in our random partial orders.  Pr\"omel, Steger and Taraz produced a series of two
papers~\cite{PST1,PST2}: in~\cite{PST1}, they demonstrated that, whatever constraint is imposed, the typical partial order arising is a $k$-layer partial
order for some~$k$.  In~\cite{PST2}, they studied the optimisation problem of finding the optimal structure given the target proportion of relations.
It follows that, if we give each partial order $P$ a weight
proportional to $e^{-cr(P)}$, where $c$ is a positive constant and 
$r(P)$ is the number of comparable pairs, then the resulting partial order will still almost surely be layered.

One way to assign weights in such a way that the majority of the weight is not supported on the class of layered posets is to give most posets
weight zero, i.e., to restrict attention to some much smaller class of partial orders.
Specifically, suppose we concentrate on the class of partial orders that can be embedded in $\M^d$, i.e., the $(d-1)$-sphere orders, and then consider
the uniform measure on this class.  What does a uniformly random $n$-element $(d-1)$-sphere order look like?  By definition, such a partial order would embed
into $\M^d$, but there is no a priori reason to suppose that it can be embedded faithfully into any finite region.  However, one reason
for optimism is that such a result is true in the case $d=2$.

As we saw in Section~\ref{sec:two}, the partial orders that can be embedded in $\M^2$ are exactly the 2-dimensional orders, and the model of random
partial orders given by a Poisson process of intensity $n$ in $\M^2$ is equivalent to that given by taking two uniformly random linear orders on $n$
labelled points and intersecting them to form a random 2-dimensional partial order.  The difference between choosing two linear orders uniformly and choosing
a uniform random 2-dimensional partial order is that some 2-dimensional partial orders have more than one representation as the
intersection of two linear orders.  The extreme case is the antichain, which has $n!$ representations as the intersection of one linear order and its reverse.  However, it turns out that most 2-dimensional partial
orders, in either model, have rather few representations.  Given a pair of linear orders on a fixed $n$-element ground-set, we say that a pair of elements
is {\em swappable} if they appear consecutively in both linear orders, and they are incomparable in the intersection.  A swappable pair
can clearly be swapped in each order, and the intersection will be unaffected.  Similarly, if a pair of linear orders has $k$ disjoint swappable pairs,
then there are (at least) $2^k-1$ other pairs of linear orders whose intersection is the same as the original pair.  If we choose the two linear orders
uniformly at random, the expected number of swappable pairs is asymptotic to~1 as $n \to \infty$, and it is easy to see that the limiting distribution is
Poisson with mean~1.  A result of Winkler~\cite{Winkler2d} states that, up to swappable pairs, most 2-dimensional partial orders have a
unique representation.  This implies that the two models of random partial orders are {\em contiguous}, i.e., any event with positive probability in one model
also has positive probability in the other.

The situation is even nicer for {\em unlabelled} 2-dimensional partial orders.  The following result was proved by El-Zahar and Sauer~\cite{EZS} and
stated in this form by Winkler~\cite{Winkler2d}.  We consider two different probability measures on the space $D_n$ of $n$-element unlabelled 2-dimensional
partial orders; one is obtained by intersecting two random linear orders on $\{1,\dots, n\}$, and the other is the
uniform probability measure on $D_n$.  In the unlabelled case, in either model, almost every element of $D_n$ has a unique representation as an intersection
of two linear orders, and, if $\Phi$ is any isomorphism-invariant statement about partial orders with a limiting probability in either model, a limiting
probability exists in the other case as well and the two limits are equal.

The implications of this result for causal sets were explored by Brightwell, Henson and Surya~\cite{BHS}.  The result of El-Zahar and Sauer
implies that a 2-dimensional partial order chosen uniformly at random is (effectively) distributed as the partial order induced on a Poisson process in
$[0,1]^2$ (or, equivalently, an interval in $\M^2$), and so such a random partial order can be embedded faithfully in an interval in $\M^2$.
This indicates how an approximation of a continuum structure can arise from a probability measure -- the uniform measure on $D_n$ --
whose definition does not refer to any continuous structure.

For each $d > 2$, it is unknown whether the model of a random $d$-dimensional partial order obtained by intersecting $d$ uniformly random linear orders
is contiguous with the model obtained by choosing a $d$-dimensional partial order uniformly at random.  It is also unknown whether the model of taking a
uniformly random order embeddable in $\M^d$ is contiguous with a Poisson process in any subset of $\M^d$.

Just to indicate some of the issues arising, let $R$ be any region of $\M^2$, and consider a partial order chosen uniformly at random from those embeddable
in~$R$.  Whatever $R$ we choose, these are exactly the 2-dimensional partial orders; we have seen that a uniformly chosen 2-dimensional partial order can
be ``faithfully'' embedded in $[0,1]^d$, and therefore in any rectangular region of the space, but not in a region of any other shape.

\end{document}